\documentclass[11pt,a4paper]{amsart}
\usepackage[T1]{fontenc}
\usepackage[utf8]{inputenc}
\usepackage{lmodern}
\usepackage{microtype}
\microtypesetup{expansion=false}
\usepackage{amsmath,amssymb,amsthm,mathtools}
\usepackage{enumitem}
\usepackage{xcolor}
\usepackage{hyperref}
\hypersetup{colorlinks=true,linkcolor=blue!45!black,citecolor=blue!45!black,urlcolor=blue!45!black,hypertexnames=false,pdftitle={Reconstructing Numerical Semigroups from Circulant Spectral Dimensions},pdfauthor={Juan Diego Velez Caicedo and Carlos A. Cadavid}}
\allowdisplaybreaks[2]
\newtheorem{theorem}{Theorem}[section]
\newtheorem{lemma}[theorem]{Lemma}
\newtheorem{proposition}[theorem]{Proposition}
\newtheorem{corollary}[theorem]{Corollary}
\theoremstyle{definition}
\newtheorem{example}[theorem]{Example}
\newtheorem{remark}[theorem]{Remark}
\newcommand{\Q}{\mathbb Q}
\newcommand{\Z}{\mathbb Z}

\newcommand{\glob}{\mathrm{glob}}
\newcommand{\ASg}{A_S^{\glob}}

\newcommand{\one}{\mathbf 1}
\newcommand{\gr}{\operatorname{gr}}
\newcommand{\Cay}{\operatorname{Cay}}

\newcommand{\rank}{\operatorname{rank}}

\title[Reconstructing Numerical Semigroups]{Reconstructing Numerical Semigroups from Circulant Spectral Dimensions}
\author{Juan Diego V\'elez Caicedo}
\address{Departamento de Matem\'aticas, Universidad Nacional de Colombia, Sede Medell\'in, Colombia}
\email{jdvelez@unal.edu.co}
\author{Carlos A. Cadavid}
\address{Departamento de Ciencias Matem\'aticas, Universidad EAFIT, Medell\'in, Colombia}
\email{ccadavid@eafit.edu.co}
\date{}

\begin{document}

\begin{abstract}
The polynomial algebra generated by the spectral moments of a family
of circulant graphs consists of all eventually polynomial sequences
in the graph order; as an unfiltered algebra, it therefore carries
little information about the jumps. We show that its natural
power-sum weight filtration retains substantially finer arithmetic
data. For a finite jump set $S$ with $\gcd(S)=1$ and $M=\max S$,
the relevant invariant is the numerical semigroup
$\Sigma_S=\langle M,M-s,M+s:s\in S\rangle$, which arises from
the lengths of closed walks winding around the cycle.
We prove that the eventual global filtered dimension is
$MD+\lfloor D/2\rfloor+1-M(M+1)-g(\Sigma_S)$.
Restricting graph orders to $n\equiv a\pmod q$, for $q\mid M$,
replaces the genus term by the number of gaps of $\Sigma_S$ in
the class $-a\pmod q$; the class $a=0$ gives the genus of the
quotient semigroup $\Sigma_S/q$. These formulas yield an inverse
reconstruction theorem: the unrestricted dimensions determine
$M$, and the dimensions in all classes modulo $M$ at one common
stable degree determine the Ap\'ery set and hence $\Sigma_S$.
We give effective thresholds, exact examples, and a finite-rank
algorithm valid in every degree. Reconstruction of the semigroup
need not identify the original jump set.
\end{abstract}
\maketitle

\section{Introduction: an inverse problem for circulant spectra}
The spectrum of a circulant graph is explicitly accessible through
Fourier analysis, and its moments count closed walks. For a family of
such graphs with a fixed jump set and varying order, the moment
functions contain an additional feature: a walk may close on the
cycle even though its displacement on the integer line is nonzero.
These \emph{winding walks} occur at exceptional graph orders, and
their first appearance depends on the arithmetic of the jumps.

Our aim is to understand how much of that arithmetic survives a
considerable loss of information. We retain only the dimensions of
spaces of polynomial expressions in the spectral moments, subject
to a bound on their natural weight. The resulting inverse problem
is concrete:
\begin{quote}
\emph{What can be reconstructed about the jumps from the filtered
spectral dimensions, and what additional information becomes
visible when the graph order is restricted to an arithmetic
progression?}
\end{quote}

Let $S$ be a nonempty finite subset of $\mathbb Z_{>0}$, and put
$r=|S|$ and $M=\max S$. For each integer $n>2M$, the undirected
circulant graph
\[
\Gamma_{n,S}=\Cay\bigl(\mathbb Z/n\mathbb Z,\{\pm s:s\in S\}\bigr)
\]
has adjacency matrix $A_{n,S}$. Its eigenvalues can be written as
$2L_S(\zeta_n^k)$, where
\[
L_S(z)=\frac12\sum_{s\in S}(z^s+z^{-s}),\qquad
\zeta_n=e^{2\pi i/n},\qquad 0\le k<n.
\]
Consequently,
\[
2^{-h}\operatorname{tr}(A_{n,S}^h)
=\sum_{k=0}^{n-1}L_S(\zeta_n^k)^h.
\]
Here $\operatorname{tr}(A_{n,S}^h)$ counts closed walks of
length $h$ with a distinguished starting vertex; the factor
$2^{-h}$ is the normalization coming from $L_S$.

We remove the trivial Fourier mode and set
\begin{equation}\label{eq:G}
G_h(n)=\sum_{k=1}^{n-1}L_S(\zeta_n^k)^h
=2^{-h}\operatorname{tr}(A_{n,S}^h)-r^h.
\end{equation}
Let $\Lambda_{\mathbb Q}=\mathbb Q[p_1,p_2,\ldots]$ be the
symmetric-function ring in its power-sum presentation, with
$\deg p_h=h$. The evaluation homomorphism
\[
\operatorname{ev}_S:\Lambda_{\mathbb Q}
\longrightarrow\mathbb Q^{\{n>2M\}},\qquad p_h\longmapsto G_h
\]
gives the filtered algebra
\[
A_S^{\glob}=\mathbb Q[G_1,G_2,\ldots]
\subset\mathbb Q^{\{n>2M\}},\qquad
\rho_S^{\glob}(D)=\dim_{\mathbb Q}A_{S,\le D}^{\glob},
\]
where $A_{S,\le D}^{\glob}$ is the image of weight at most $D$.
The weight assigns cost $h$ to the basic statistic $G_h$,
reflecting its walk length; products have the sum of the
corresponding costs.

The use of a filtration is essential. In fact, the unfiltered
image $A_S^{\glob}$ consists of \emph{all} rational-valued
sequences on the admissible graph orders that eventually agree
with a polynomial in $n$
(Proposition~\ref{prop:glob-image}). Thus, apart from the
choice of its index set, this algebra does not distinguish jump
sets. The dimensions $\rho_S^{\glob}(D)$ measure something
finer: the cost, in moment weight, of producing the required
polynomial and its exceptional finite corrections.

The mechanism behind this dependence is elementary but revealing.
By Fourier orthogonality, $G_h(n)$ is affine in $n$ for $n>hM$.
At smaller orders the deviation is a finite correction supported
where a signed $S$-walk has displacement $kn$ for some nonzero
integer $k$. Denote by $\nu_S(n)$ the minimum length of such a
walk. A walk with displacement $n$ and length $h$ satisfies
\[
hM-n=\sum_{i=1}^h(M-\varepsilon_i s_i),
\qquad \varepsilon_i\in\{-1,1\},\quad s_i\in S.
\]
The possible step deficits are therefore $M-s$ and $M+s$.
Adding the generator $M$ keeps track of the unused steps when
one considers an upper bound on the winding length. This leads to
the numerical semigroup
\[
\Sigma_S=\langle M,M-s,M+s:s\in S\rangle.
\]
The zero element arising from $M-s$ when $s=M$ can be
omitted from the displayed generating set.
When $\gcd(S)=1$, its complement in $\mathbb Z_{\ge0}$ is finite.
Write $g_S$ for this genus and $\delta_t$ for the least member of
$\Sigma_S$ congruent to $t$ modulo $M$. The representatives
$\delta_0,\ldots,\delta_{M-1}$ form its Ap\'ery set.
We prove that, for all sufficiently large $n$,
\[
\nu_S(n)=\frac{n+\delta_{-n\bmod M}}{M}.
\]
This identity explains why the gaps of $\Sigma_S$ appear in
spectral dimensions: at a fixed walk length, each gap removes
one possible exceptional graph order. The weight filtration
measures precisely the resulting loss of independent corrections.

The first theorem makes this connection quantitative.
\begin{theorem}[Eventual full-spectrum dimension]\label{thm:growth}
If $\gcd(S)=1$, then, for every sufficiently large $D$,
\begin{equation}\label{eq:global-asymptotic}
\rho_S^{\glob}(D)
=MD+\left\lfloor\frac D2\right\rfloor+1-M(M+1)-g_S.
\end{equation}
An explicit effective stable-range bound is provided in
Theorem~\ref{thm:global-exact}.
\end{theorem}

This formula contains an instructive limitation: once $M$ is fixed,
the eventual unrestricted dimension depends on $S$ through the
single integer $g_S$. Distinct semigroups of the same genus are
therefore indistinguishable by these dimensions. To recover the
individual gaps, we keep track of the graph orders at which the
corrections occur. For $q\mid M$ and
$0\le a<q$, let $\rho_{S;q,a}(D)$ be the dimension of the image
of $A_{S,\le D}^{\glob}$ after restriction to graph orders
$n>2M$ satisfying $n\equiv a\pmod q$. Write
\[
g_b(T;q)=\#\{v\in\mathbb Z_{\ge0}\setminus T:v\equiv b\pmod q\},
\qquad T/q=\{j\ge0:qj\in T\}.
\]
\begin{theorem}[Congruence-restricted spectral reconstruction]
\label{thm:residue-intro}
Assume $\gcd(S)=1$. For each positive divisor $q\mid M$,
each $0\le a<q$, and all sufficiently large $D$,
\begin{equation}\label{eq:residue-intro}
\rho_{S;q,a}(D)
=\frac Mq D+\left\lfloor\frac D2\right\rfloor+1
-\frac Mq\left(\frac Mq+1\right)-g_{-a\bmod q}(\Sigma_S;q).
\end{equation}
The stable-range threshold may be chosen effectively and uniformly
in $q$ and $a$. In particular the constant defect at orders divisible
by $q$ is $g(\Sigma_S/q)$. The $M$ restricted dimensions modulo $M$
at one common stable degree determine the full Ap\'ery set of
$\Sigma_S$, and therefore the semigroup itself.
\end{theorem}

The inverse consequence is especially explicit. Because the
unrestricted dimension grows by $2M+1$ when the weight increases
by two, we can first recover the previously unknown maximum jump.
Next, taking $q=M$ isolates one gap count per residue class:
\begin{theorem}[Inverse reconstruction]
\label{thm:inverse-intro}
Assume $\gcd(S)=1$. For every sufficiently large integer $D$,
the unrestricted and congruence-restricted dimensions determine
the maximum jump $M$ and every Ap\'ery representative $\delta_t$:
\begin{equation}\label{eq:inverse-intro}
\boxed{\begin{aligned}
M&=\frac{\rho_S^{\glob}(D+2)-\rho_S^{\glob}(D)-1}{2},\\[3pt]
\delta_t&=t+M\left(D+\left\lfloor\frac D2\right\rfloor-1
-\rho_{S;M,-t\bmod M}(D)\right),\quad 0\le t<M.
\end{aligned}}
\end{equation}
In particular, since
$\Sigma_S=\bigcup_{t=0}^{M-1}(\delta_t+M\mathbb Z_{\ge0})$,
these dimensions reconstruct the \emph{defect semigroup}
$\Sigma_S$. In general they do not determine the original
jump set $S$.
\end{theorem}
The proof and an explicit example of nonuniqueness are given in
Corollary~\ref{cor:inverse-reconstruction}.

Here is a concrete reason to retain the congruence restrictions.
Take $S_1=\{1,2,6\}$ and $S_2=\{1,3,6\}$. Their defect
semigroups have different gap sets, $\{1,2,3\}$ and $\{1,2,4\}$,
but the same genus $3$. Therefore the eventual unrestricted
dimension formulas are identical. Restriction to even graph orders
separates the two cases: the eventual dimensions differ by one.
We give the exact calculation in Example~\ref{ex:equal-genus}.

The coprimality hypothesis is needed for the finite-gap
description. For completeness, if $a=\gcd(S)$ and $M_0=M/a$,
then the unrestricted dimension still has a linear asymptotic:
\begin{equation}\label{eq:global-noncoprime}
\rho_S^{\glob}(D)=
\left(M_0\sum_{d\mid a}\frac{\varphi(d)}d+\frac12\right)D+O_S(1),
\end{equation}
where $\varphi$ is Euler's totient function. This is proved in
Proposition~\ref{prop:noncoprime}. It is separate from the
inverse theorem, which assumes $\gcd(S)=1$.

The ingredients from numerical-semigroup theory are classical.
We use the theory of Ap\'ery sets as developed in~\cite{Rosales}.
Adeniran et al.~\cite{Adeniran2019} established a root-of-unity
formula for the genus of a quotient $T/q$, and Lee and
Nam~\cite{LeeNam2024} computed quotient genera for several
families, including semigroups generated by arithmetic
progressions. Our earlier cyclotomic evaluation framework
\cite{VelezCadavid2026} also explains eventual polynomiality
of bounded-degree symmetric evaluations. These antecedents
motivate the present question but do not, by themselves,
provide the dimensions of the spectral evaluation images
in each weight.

The contributions of this paper are threefold. First, we
construct filtration-compatible bases for the full-spectrum
and restricted evaluation spaces, with effective stable-range
bounds. Second, we identify the exact spectral deficits with
gap counts, including quotient genera, and invert the formulas
to reconstruct the Ap\'ery set. Third, we prove all-degree
rank formulas in representative families and give a finite
evaluation-matrix procedure that also applies outside the
coprime and divisibility hypotheses. The inverse problem
therefore has a precise answer at the level of the defect
semigroup, together with a precise limitation on recovering
the jump set itself.

Section~2 derives the moment and winding identities. Section~3 proves
the dimension and reconstruction theorems, including effective
thresholds and quotient-genus formulas. Section~4 contains exact
examples and a finite-rank computation method.

\section{Fourier moments and winding corrections}
Expand the Laurent polynomial
\begin{equation}\label{eq:Laurent}
L_S(z)^h=\sum_{|m|\le hM}b_{h,m}z^m,\qquad a_h=b_{h,0}.
\tag{2.1}\end{equation}
Every coefficient $b_{h,m}$ is a nonnegative rational number,
$b_{h,m}=b_{h,-m}$, and
\[
a_1=0,\quad a_2=r/2>0,\quad
\sum_m b_{h,m}=r^h.
\]
\begin{lemma}[Full-spectrum moment formula]\label{lem:global-bridge}
For every $n>2M$ and $h\ge1$,
\begin{equation}\label{eq:G-bridge}
G_h(n)=a_hn-r^h+n\sum_{\substack{m\ne0\\n\mid m}}b_{h,m}.
\tag{2.2}\end{equation}
In particular, $G_h(n)=a_hn-r^h$ when $n>hM$.
\end{lemma}
\begin{proof}
Interchanging two finite sums in~\eqref{eq:G} yields
\[
\sum_{k=0}^{n-1}L_S(\zeta_n^k)^h
=\sum_m b_{h,m}\sum_{k=0}^{n-1}\zeta_n^{km}
=n\sum_{n\mid m}b_{h,m}.
\]
The omitted trivial Fourier mode contributes $L_S(1)^h=r^h$.
Splitting off the $m=0$ term proves the formula.
\end{proof}

The coefficient $b_{h,m}$ is the normalized number of signed
$S$-walks of length $h$ and total displacement $m$. Thus the
finite correction in~\eqref{eq:G-bridge} records precisely those
walks winding a nonzero number of times around the cycle.

\section{The global algebra and its numerical semigroup}
The full-spectrum moments have a simple polynomial part. It is convenient to separate it from the terms coming from winding walks. Put
\begin{equation}\label{eq:tau-Q}
\tau_S=\frac2r(G_2+r^2)=n,\qquad
Q_h(n)=G_h(n)-a_hn+r^h
=n\sum_{\substack{m\ne0\\n\mid m}}b_{h,m},\quad h\ge3.
\tag{3.1}\end{equation}
Here $\tau_S=n$ means equality as arithmetic functions on the admissible levels. The correction $Q_h$ is nonnegative and supported in the finite interval $2M<n\le hM$.

\begin{proposition}[The complete global image]\label{prop:glob-image}
For every nonempty finite $S$, the algebra $\ASg$ is exactly the algebra of rational-valued sequences on $n>2M$ that are eventually polynomial in $n$. An element of filtered weight at most $D$ agrees, for all sufficiently large $n$, with a polynomial of degree at most $\lfloor D/2\rfloor$.
\end{proposition}
\begin{proof}
For $n>DM$, each moment $G_h$ with $h\le D$ has the form $a_hn-r^h$. Since $G_1$ has zero slope, every monomial of total weight at most $D$ gives a polynomial of degree at most $\lfloor D/2\rfloor$. Conversely, the function $\tau_S=n$ belongs to our algebra, so all polynomial sequences belong to it.

We also have all finitely supported sequences. Indeed, for $m>2M$ the coefficient $b_{m,mM}=2^{-m}$ is positive, and hence $Q_m(m)>0$. As $Q_m$ has finite support, interpolation gives the indicator of the level $m$:
\[
E_{m,S}(n)=\frac{Q_m(n)}{Q_m(m)}
\prod_{\substack{2M<j\le mM\\j\ne m}}\frac{n-j}{m-j}.
\]
This belongs to $\ASg$. Combining finitely many such indicators with a polynomial produces any eventually polynomial sequence.
\end{proof}

We now determine how much weight is needed to obtain these finite corrections. Define the first nonzero winding length by
\[
\nu_S(n)=\min\left\{h\ge1:\sum_{i=1}^h\varepsilon_i s_i\in n\Z\setminus\{0\},\quad
\varepsilon_i\in\{\pm1\},\ s_i\in S\right\}.
\]
Since $n>2M$, we have $\nu_S(n)\ge3$. The coefficients of $L_S^h$ are nonnegative, and therefore
\[
\nu_S(n)=\min\{h\ge3:Q_h(n)>0\}.
\]
Write $X_h=\{n>2M:\nu_S(n)\le h\}$, $B_h=X_h\setminus X_{h-1}$ and $b_h=|B_h|$, with $X_2=\varnothing$.

\begin{lemma}[Triangular interpolation basis]\label{lem:triangular-global}
For each $H\ge3$, the functions
\begin{equation}\label{eq:triangular-global}
\{n^jQ_h(n):3\le h\le H,\quad 0\le j<b_h\}
\tag{3.2}\end{equation}
form a $\Q$-basis of all rational-valued functions supported on $X_H$. The function $n^jQ_h$ can be obtained in weight at most $h+2j$.
\end{lemma}
\begin{proof}
If $n\in B_h$, then $\nu_S(n)=h$. Consequently $Q_h(n)>0$,
whereas $Q_j(n)=0$ for every $j<h$. On the $b_h$ distinct
points of $B_h$, the evaluation matrix of
$Q_h,nQ_h,\ldots,n^{b_h-1}Q_h$ is the product of the
Vandermonde matrix $(n^j)$ and the invertible diagonal matrix
with diagonal entries $Q_h(n)$. Its determinant is nonzero.

To interpolate an arbitrary function on $X_H$, first match its
values on $B_H$ with a unique linear combination of these $b_H$
functions, and subtract the combination. The remaining function
vanishes on $B_H$. Repeat on $B_{H-1}$ and proceed down to $B_3$.
At each step, the newly used functions vanish on all higher layers,
so previously matched values are unchanged. This proves spanning.
If a linear combination of the displayed family is zero, restriction
first to $B_H$ forces all coefficients belonging to $h=H$ to vanish;
descending through the layers proves independence.

Finally, $n=\tau_S$ has weight two, and
$Q_h=G_h-a_h\tau_S+r^h$ has weight at most $h$ for $h\ge3$.
Thus $n^jQ_h$ belongs to weight at most $h+2j$.
\end{proof}

From now until Theorem~\ref{thm:global-exact}, suppose that $\gcd(S)=1$. Consider
\begin{equation}\label{eq:semigroup}
\Sigma_S=\langle M,M-s,M+s:s\in S\rangle.
\tag{3.3}\end{equation}
This is a numerical semigroup. For $0\le t<M$, let
\[
\delta_t=\min\{\delta\in\Sigma_S:\delta\equiv t\pmod M\},\qquad
g_S=\sum_{t=0}^{M-1}\left\lfloor\frac{\delta_t}{M}\right\rfloor.
\]
The $\delta_t$ form the Ap\'ery set of $\Sigma_S$ with respect to $M$; we refer to~\cite{Rosales} for the standard facts about numerical semigroups.

For completeness, we specify an effective bound. Choose for each $t$ a representation
\[
\delta_t=\sum_{i=1}^{\ell_t}(M-\varepsilon_i s_i),\qquad
\varepsilon_i\in\{\pm1\},\quad s_i\in S,
\]
with no summand equal to $M$, and put $\ell_0=0$. Set
\[
\begin{aligned}
\delta^*&=\max_t\delta_t,\qquad \ell^*=\max_t\ell_t,\\
N^*&=\max\{2M,\delta^*,M\ell^*\},\qquad H=\max\{3,N^*\}.
\end{aligned}
\]
The sets $B_h$ and the integers $b_h$ for $3\le h\le H$ are determined by finite signed-walk calculations. Define
\[
D_0=\max_{\substack{3\le h\le H\\b_h>0}}\{h+2(b_h-1)\},
\qquad D^*=\max\{D_0,H+2M-2\}.
\]
The set indexing $D_0$ is nonempty because $3M\in B_3$.

\begin{lemma}[Ap\'ery formula for winding]\label{lem:apery-winding}
For $n>N^*$, take $0\le t<M$ with $t\equiv-n\pmod M$. Then
\begin{equation}\label{eq:winding-formula}
\nu_S(n)=\frac{n+\delta_t}{M}.
\tag{3.4}\end{equation}
For every $h\ge H$,
\begin{equation}\label{eq:X-cardinality}
|X_h|=Mh-2M-g_S,
\tag{3.5}\end{equation}
and for $h>H$ the $M$ elements of $B_h$ are exactly
$hM-\delta_t$, $0\le t<M$.
\end{lemma}
\begin{proof}
Fix $n>N^*$ and put $t\equiv-n\pmod M$, with $0\le t<M$.
The integer
\[
h_0=\frac{n+\delta_t}{M}
\]
is positive. We prove first that every signed walk with nonzero
winding has length at least $h_0$.

Suppose a signed walk of length $h$ has displacement $kn$, where
$k\in\mathbb Z\setminus\{0\}$. By reversing all signs, we may assume
$k\ge1$. If $k=1$, its signed steps $\varepsilon_i s_i$ satisfy
\[
Mh-n=\sum_{i=1}^h(M-\varepsilon_i s_i)\in\Sigma_S.
\]
Since $Mh-n\equiv t\pmod M$, minimality of $\delta_t$ yields
$Mh-n\ge\delta_t$, or $h\ge h_0$. If $k\ge2$, then
$hM\ge kn\ge2n$. But $n>N^*\ge\delta^*\ge\delta_t$,
so $2n>n+\delta_t=Mh_0$, and hence $h>h_0$.

Conversely, use the chosen representation
\[
\delta_t=\sum_{i=1}^{\ell_t}(M-\varepsilon_i s_i).
\]
An Ap\'ery representative has a representation without a summand
equal to $M$: deleting such a summand would give a smaller element
in the same residue class. The chosen signed steps have total
displacement $\ell_tM-\delta_t$. Since
$n>N^*\ge M\ell^*$, we have $h_0\ge\ell_t$. Append
$h_0-\ell_t$ steps of size $+M$. The resulting walk has length
$h_0$ and displacement
\[
(\ell_tM-\delta_t)+(h_0-\ell_t)M=h_0M-\delta_t=n.
\]
Thus $\nu_S(n)=h_0$, proving~\eqref{eq:winding-formula}.

We now establish the support identity
\begin{equation}\label{eq:stable-support-identity}
X_h=\{n>2M:Mh-n\in\Sigma_S\},\qquad h\ge H.
\end{equation}
For $n>N^*$, this follows from~\eqref{eq:winding-formula}
and the description of the members of $\Sigma_S$ in residue class
$t$ as $\delta_t+M\mathbb Z_{\ge0}$. If $2M<n\le N^*$,
then $\nu_S(n)\le n\le H$, since $n$ consecutive $+M$
steps have displacement $nM\in n\mathbb Z\setminus\{0\}$.
Also $n+\delta_t\le2N^*\le MH$ when $M\ge2$;
for $M=1$ we have $\delta_0=0$ and the same inequality is
immediate. Thus $Mh-n\ge\delta_t$ for $h\ge H$, and its
residue is $t$, so $Mh-n\in\Sigma_S$. This proves the identity.

Under $v=Mh-n$, the right-hand side of
\eqref{eq:stable-support-identity} is in bijection with
\[
\{v\in\Sigma_S:0\le v<M(h-2)\}.
\]
The largest gap of $\Sigma_S$, when gaps exist, is
$\max_t\delta_t-M=\delta^*-M$. Our choice of $H$ ensures that
every such gap is below $M(h-2)$ for $h\ge H$.
Consequently
\[
|X_h|=M(h-2)-g_S= Mh-2M-g_S.
\]
Taking successive differences gives $|B_h|=M$ for $h>H$.
For each $t$, the order $n=hM-\delta_t$ is admissible and
larger than $N^*$ when $h>H$; by~\eqref{eq:winding-formula}
its first winding length is $h$. These $M$ orders lie in
distinct residue classes modulo $M$, so they constitute $B_h$.
\end{proof}

\begin{theorem}[Exact eventual global filtered dimension]\label{thm:global-exact}
If $\gcd(S)=1$, formula~\eqref{eq:global-asymptotic} holds for every $D\ge D^*$. More precisely, for every $D\ge D_0$ a basis of $A_{S,\le D}^{\glob}$ is
\begin{equation}\label{eq:global-basis}
\begin{aligned}
&\{1,n,\ldots,n^{\lfloor D/2\rfloor}\}\\[-2pt]
&\quad\cup\{n^jQ_h:3\le h\le H,\ 0\le j<b_h\}\\[-2pt]
&\quad\cup\{n^jQ_h:h>H,\ 0\le j<M,\ h+2j\le D\}.
\end{aligned}
\tag{3.6}\end{equation}
In particular, the Hilbert series of $\gr\ASg$ is rational, with denominator dividing $1-z^2$.
\end{theorem}
\begin{proof}
Every displayed function belongs to the indicated filtration:
$n$ has weight two, and $Q_h$ has weight at most $h$.
The finite middle family is therefore available for $D\ge D_0$.
Lemma~\ref{lem:triangular-global} gives independence of all correction
functions, because their restrictions are triangular on the layers
$B_h$. A polynomial in $n$ cannot equal a nonzero finitely supported
function: restricting to an infinite tail forces the polynomial to
vanish identically. Thus the entire proposed family is independent.

We prove spanning with explicit control of the weight. Expand every
$G_h$ as $a_hn-r^h+Q_h$, and expand a product of moments of total
weight at most $D$. Its term containing no $Q_h$ is a polynomial
in $n$ of degree at most $\lfloor D/2\rfloor$, since $a_1=0$ and
each occurrence of $n$ therefore comes from a moment of weight
at least two.

Every other term has the form
\[
c\,n^j Q_{h_1}\cdots Q_{h_m},\qquad
m\ge1,\quad h_i\ge3,
\]
with $2j+\sum_i h_i\le D$. Put $h=\min_i h_i$. The support of
this term is contained in $X_h$. If $h\le H$, the triangular
interpolation lemma expresses it in the middle family, whose
weights are at most $D_0\le D$.

Suppose $h>H$. The same interpolation procedure on $X_h$ uses
only the middle family and the functions
\[
n^{j'}Q_k,\qquad H<k\le h,\quad 0\le j'<M.
\]
All the latter have weight at most $h+2(M-1)$.
If $m=1$ and $j<M$, the original term $n^jQ_h$ is already
listed in~\eqref{eq:global-basis}. If $m=1$ and $j\ge M$,
then $h+2(M-1)\le h+2j\le D$. Finally, if $m\ge2$,
then $D\ge\sum_i h_i\ge2h$, and
$h+2(M-1)\le2h$ because $h>H\ge2M$.
Together with $D\ge D_0$, these inequalities show that
interpolation never leaves weight at most $D$. This proves spanning.

For $D\ge D^*$ the low correction family has $|X_H|=MH-2M-g_S$ elements. The high correction family contributes
\[
\sum_{j=0}^{M-1}(D-H-2j)=M(D-H)-M(M-1)
\]
elements, and the polynomial family has $\lfloor D/2\rfloor+1$ elements. Adding the three numbers gives~\eqref{eq:global-asymptotic}. Eventually, the successive graded dimensions alternate between $M+1$ in even degrees and $M$ in odd degrees; the rationality assertion follows.
\end{proof}

\begin{proposition}[The noncoprime leading coefficient]
\label{prop:noncoprime}
For every nonempty finite $S$, without assuming $\gcd(S)=1$,
the asymptotic formula~\eqref{eq:global-noncoprime} holds.
\end{proposition}
\begin{proof}
Put $a=\gcd(S)$, $S'=S/a$, and $M_0=M/a$. Let
$g=\gcd(n,a)$ and $m=n/g$. For any integer displacement $v$,
\[
n\mid av\quad\Longleftrightarrow\quad m\mid v.
\]
Thus $\nu_S(n)=\nu_{S'}(m)$. Conversely, a positive integer
$n$ has $\gcd(n,a)=g$ and reduced modulus $m$ precisely when
$n=gm$, $g\mid a$, and $\gcd(m,a/g)=1$.
The admissibility condition $n>2M$ removes only a bounded initial
set of reduced moduli for each fixed divisor $g$.

Since $\gcd(S')=1$, Lemma~\ref{lem:apery-winding} shows that,
for large $h$, the set $X_h(S')$ differs from the interval
$2M_0<m\le hM_0$ by a bounded number of elements. For every
fixed positive integer $d$,
\[
\#\{1\le m\le x:\gcd(m,d)=1\}
=\frac{\varphi(d)}d x+O_d(1),
\]
as follows directly by inclusion--exclusion over the divisors of
$d$. Applying this count to $d=a/g$ and summing over $g\mid a$
gives
\[
|X_h(S)|
=M_0\left(\sum_{d\mid a}\frac{\varphi(d)}d\right)h
+O_S(1).
\]
Hence $b_h=|X_h|-|X_{h-1}|$ is bounded independently of $h$.
Choose an integer $B$ with $b_h\le B$ for every $h\ge3$,
and set $C=2(B-1)$. By Lemma~\ref{lem:triangular-global},
every rational-valued function supported on $X_{D-C}$ belongs to
the weight-$D$ image for sufficiently large $D$, because each of
its basis terms $n^jQ_h$ has weight at most
$h+2(b_h-1)\le D$. Conversely, an expression of weight at most
$D$, after subtracting its eventual polynomial part, is supported
on $X_D$. Its polynomial part has degree at most
$\lfloor D/2\rfloor$.
The polynomial and finite-support subspaces intersect trivially,
so
\[
\left\lfloor\frac D2\right\rfloor+1+|X_{D-C}|
\le\rho_S^{\glob}(D)\le
\left\lfloor\frac D2\right\rfloor+1+|X_D|.
\]
The preceding support estimate now proves
\eqref{eq:global-noncoprime}.
\end{proof}

\subsection{Congruence restrictions and quotient genera}\label{subsec:quotients}
For the rest of this section, assume again that $\gcd(S)=1$, and retain the effective integers $H,D_0,D^*$ defined above. Let $q$ be a positive divisor of $M$, let $0\le a<q$, and put
\[
I_{q,a}=\{n>2M:n\equiv a\pmod q\},\qquad N=M/q.
\]
We give the restriction image its induced filtration:
\[
A_{S;q,a,\le D}^{\glob}
=\operatorname{res}_{I_{q,a}}(A_{S,\le D}^{\glob}),
\qquad \rho_{S;q,a}(D)=\dim_\Q A_{S;q,a,\le D}^{\glob}.
\]
Only the graph orders are restricted; the moment weights stay unchanged. We use the gap counts $g_b(T;q)$ defined in the introduction, with residue subscripts interpreted modulo $q$.

\begin{lemma}[Winding counts in a residue class]\label{lem:residue-counts}
For every $h\ge H$,
\begin{equation}\label{eq:residue-counts}
|X_h\cap I_{q,a}|=Nh-2N-g_{-a}(\Sigma_S;q).
\tag{3.7}\end{equation}
Consequently $|B_h\cap I_{q,a}|=N$ for $h>H$.
\end{lemma}
\begin{proof}
The proof of Lemma~\ref{lem:apery-winding} gives, for every $h\ge H$,
\[
X_h=\{n>2M:Mh-n\in\Sigma_S\}.
\]
Under the substitution $v=Mh-n$, the restricted support is therefore described by
\[
0\le v<M(h-2),\qquad v\in\Sigma_S,\qquad v\equiv-a\pmod q.
\]
Because $q\mid M$, this interval contains exactly $N(h-2)$ nonnegative integers in each residue class modulo $q$. It also contains every gap of $\Sigma_S$: the largest gap, when one exists, is $\delta^*-M$, and our choice of $H$ gives $M(H-2)-1\ge\delta^*-M$. Subtracting the gaps in class $-a$ proves~\eqref{eq:residue-counts}. Taking consecutive differences proves the assertion about $B_h$.
\end{proof}

\begin{theorem}[Exact restricted filtered dimension]\label{thm:restricted-dimension}
For every $q\mid M$, every $0\le a<q$, and every integer $D\ge D^*$,
\begin{equation}\label{eq:restricted-dimension}
\rho_{S;q,a}(D)
=ND+\left\lfloor\frac D2\right\rfloor+1
-N(N+1)-g_{-a}(\Sigma_S;q),\qquad N=M/q.
\tag{3.8}\end{equation}
The same threshold $D^*$ works for all $q$ and $a$. The associated-graded Hilbert series of each restriction image is rational, with denominator dividing $1-z^2$.
\end{theorem}
\begin{proof}
Write $b_h'=|B_h\cap I_{q,a}|$ and suppress the restriction symbol on functions. The argument of Lemma~\ref{lem:triangular-global} applies on each restricted layer. Indeed, $Q_h(n)>0$ on $B_h$, so the functions $n^jQ_h$, $0\le j<b_h'$, give a Vandermonde matrix multiplied by an invertible diagonal matrix on $B_h\cap I_{q,a}$. Interpolating successively from the highest layer down gives a basis of all functions supported on $X_L\cap I_{q,a}$, for any $L\ge3$.

For $D\ge D^*$, a basis of the filtered restriction image is
\begin{equation}\label{eq:restricted-basis}
\begin{aligned}
&\{1,n,\ldots,n^{\lfloor D/2\rfloor}\}\\[-2pt]
&\quad\cup\{n^jQ_h:3\le h\le H,\ 0\le j<b_h'\}\\[-2pt]
&\quad\cup\{n^jQ_h:h>H,\ 0\le j<N,\ h+2j\le D\}.
\end{aligned}
\tag{3.9}\end{equation}
Every displayed function is available in degree $D$. For the middle family this follows from $b_h'\le b_h$ and $D\ge D_0$; for the last family we use Lemma~\ref{lem:residue-counts}. The correction functions are independent by triangular interpolation. A polynomial vanishing on an infinite tail of $I_{q,a}$ is zero, so the polynomial family is independent of the finite corrections.

To check spanning, express each moment in terms of $n$ and the $Q_h$. A monomial containing a correction factor is supported on $X_h$, where $h$ is its smallest correction index. If $h\le H$, it reduces to the middle family. Suppose $h>H$. A monomial $n^jQ_h$ with $j<N$ already occurs in~\eqref{eq:restricted-basis}. For $j\ge N$, interpolation on $X_h\cap I_{q,a}$ uses high-layer weights at most $h+2(N-1)\le h+2j$. If the monomial has at least two correction factors, its weight is at least $2h$, whereas $h+2(N-1)\le2h$ because $h>H\ge2M\ge2N$. The low family is already available in degree $D$. Thus no reduction leaves the required filtered piece.

By~\eqref{eq:residue-counts}, the middle family has $NH-2N-g_{-a}(\Sigma_S;q)$ elements. Since $D\ge H+2M-2\ge H+2N-2$, the high family has
\[
\sum_{j=0}^{N-1}(D-H-2j)=N(D-H)-N(N-1)
\]
elements. Adding the polynomial dimension $\lfloor D/2\rfloor+1$ proves~\eqref{eq:restricted-dimension}. Eventually the successive graded dimensions are $N+1$ in even degree and $N$ in odd degree, giving the stated rational denominator.
\end{proof}

\begin{corollary}[Quotient genus in the spectral constant]\label{cor:quotient-genus}
For $q\mid M$, put $N=M/q$. For every $D\ge D^*$,
\begin{equation}\label{eq:quotient-dimension}
\rho_{S;q,0}(D)=ND+\left\lfloor\frac D2\right\rfloor+1
-N(N+1)-g(\Sigma_S/q).
\tag{3.10}\end{equation}
Moreover,
\begin{equation}\label{eq:quotient-apery}
\begin{aligned}
\operatorname{Ap}(\Sigma_S/q,N)&=\{\delta_{qj}/q:0\le j<N\},\\
g(\Sigma_S/q)&=\sum_{j=0}^{N-1}\left\lfloor\frac{\delta_{qj}}M\right\rfloor.
\end{aligned}
\tag{3.11}\end{equation}
\end{corollary}
\begin{proof}
The map $m\mapsto qm$ identifies the gaps of $\Sigma_S/q$ with the gaps of $\Sigma_S$ divisible by $q$. Hence $g(\Sigma_S/q)=g_0(\Sigma_S;q)$, and~\eqref{eq:quotient-dimension} is the case $a=0$ of Theorem~\ref{thm:restricted-dimension}. If $m\equiv j\pmod N$, then $qm\equiv qj\pmod M$. The least such $qm$ in $\Sigma_S$ is $\delta_{qj}$, which is divisible by $q$. This proves the Ap\'ery formula, and counting the gaps below these representatives gives the genus formula.
\end{proof}

We recall the general quotient-genus identity of Adeniran et al.~\cite[Theorem~1.3]{Adeniran2019}. For a numerical semigroup $T$, write $H_T(z)=\sum_{m\in T}z^m$ and $\zeta_q=e^{2\pi i/q}$. Then
\begin{equation}\label{eq:known-quotient-genus}
g(T/q)=\frac1q\left(g(T)+\frac{q-1}{2}
-\sum_{j=1}^{q-1}H_T(\zeta_q^j)\right).
\tag{3.12}\end{equation}
These evaluations use the rational expression
\[
H_T(z)=\frac1{1-z}-\sum_{m\notin T}z^m,
\]
which is regular at every root of unity different from one. Apparent singularities in an Ap\'ery expression with denominator $1-z^M$ must therefore be canceled before evaluating at such roots. Substituting~\eqref{eq:known-quotient-genus} into~\eqref{eq:quotient-dimension} expresses the spectral constant in terms of the numerical-semigroup Hilbert series. The semigroup series and the associated-graded spectral Hilbert series are related here through the support count~\eqref{eq:residue-counts}.

\subsection{Reconstruction of the defect semigroup}\label{subsec:reconstruction}
The full collection of congruence restrictions contains more information than the eventual global dimension formula.

\begin{corollary}[Recovery of the Ap\'ery set]\label{cor:apery-reconstruction}
For every $D\ge D^*$ and $0\le t<M$,
\begin{equation}\label{eq:apery-reconstruction}
\delta_t=t+M\left(D+\left\lfloor\frac D2\right\rfloor-1
-\rho_{S;M,-t\bmod M}(D)\right).
\tag{3.13}\end{equation}
Thus the $M$ restricted dimensions modulo $M$, taken in one stable degree, determine $\Sigma_S$.
\end{corollary}
\begin{proof}
The elements of $\Sigma_S$ in residue class $t$ are exactly $\delta_t+M\Z_{\ge0}$. Consequently
\[
g_t(\Sigma_S;M)=\frac{\delta_t-t}{M}.
\]
Taking $q=M$, $N=1$, and $a\equiv-t\pmod M$ in~\eqref{eq:restricted-dimension} gives~\eqref{eq:apery-reconstruction}. The semigroup is then recovered from
\[
\Sigma_S=\bigcup_{t=0}^{M-1}(\delta_t+M\Z_{\ge0}).
\]
\end{proof}

\begin{corollary}[Inverse reconstruction without a prescribed maximum jump]
\label{cor:inverse-reconstruction}
Let $S$ be a nonempty finite set with $\gcd(S)=1$.
Its eventual global and congruence-restricted filtered dimensions
determine both $M=\max S$ and the numerical semigroup $\Sigma_S$.
More explicitly, for any $D\ge D^*$,
\begin{equation}\label{eq:recover-M}
M=\frac{\rho_S^{\glob}(D+2)-\rho_S^{\glob}(D)-1}{2}.
\end{equation}
Having recovered $M$, the $M$ dimensions
$\rho_{S;M,a}(D)$, $0\le a<M$, determine the Ap\'ery
representatives by~\eqref{eq:apery-reconstruction}.

If $S$ and $S'$ are two coprime jump sets with the same maximum $M$,
then their eventual restricted dimensions for all residue classes
modulo $M$ agree if and only if $\Sigma_S=\Sigma_{S'}$.
Neither these dimensions nor the recovered semigroup necessarily
determine the original jump set.
\end{corollary}
\begin{proof}
Formula~\eqref{eq:global-asymptotic} gives
$\rho_S^{\glob}(D+2)-\rho_S^{\glob}(D)=2M+1$
whenever $D\ge D^*$, establishing~\eqref{eq:recover-M}.
The second assertion follows from
Corollary~\ref{cor:apery-reconstruction}; the explicit
representatives determine the semigroup as a union of arithmetic
progressions. Equality of all the restricted eventual dimensions
forces equality of every $\delta_t$ and hence of the semigroups.
Conversely, equal semigroups with the same $M$ have the same
gap counts in every class, so
Theorem~\ref{thm:restricted-dimension} gives equal dimensions
for all sufficiently large $D$.

For the final qualification, consider $S=\{5,6\}$ and
$S'=\{1,5,6\}$. Both are coprime and have maximum jump six,
but $S\ne S'$, while both defect semigroups contain
$6-5=1$. Hence $\Sigma_S=\Sigma_{S'}=\mathbb Z_{\ge0}$,
and their eventual global and restricted dimensions agree.
\end{proof}

\begin{example}[Equal global dimensions, different even-order ranks]\label{ex:equal-genus}
Take $S_1=\{1,2,6\}$ and $S_2=\{1,3,6\}$. Their defect semigroups are
\[
\Sigma_{S_1}=\langle4,5,6,7\rangle,
\qquad\Sigma_{S_2}=\langle3,5,7\rangle,
\]
with gap sets $\{1,2,3\}$ and $\{1,2,4\}$, respectively. Both have genus three and $M=6$, so both eventual global dimension formulas are
\[
\rho_{S_i}^{\glob}(D)=6D+\lfloor D/2\rfloor-44.
\]
Their quotients by two are $\langle2,3\rangle$ and $\langle3,4,5\rangle$, whose genera are one and two. Corollary~\ref{cor:quotient-genus} therefore gives
\[
\rho_{S_1;2,0}(D)=3D+\lfloor D/2\rfloor-12,
\qquad
\rho_{S_2;2,0}(D)=3D+\lfloor D/2\rfloor-13
\]
for all sufficiently large $D$.
\end{example}

\begin{remark}[The divisibility hypothesis]
The condition $q\mid M$ makes the stable restricted layer size constant. Without this condition additional periodic terms may occur. For $S=\{1\}$, restrict to graph orders divisible by three. Each level $3j$ first appears in weight $3j$, and one triangular correction is needed for each such level. The restricted dimension is
\[
\rho(D)=\left\lfloor\frac D2\right\rfloor
+\left\lfloor\frac D3\right\rfloor+1.
\]
Here $\Sigma_S=\Z_{\ge0}$ and its quotient by three still has genus zero, while the dimension has a period-three contribution.
\end{remark}

\subsection{A family generated by an arithmetic progression}\label{subsec:ap-family}
For two coprime jumps, the defect semigroup is generated by a three-term arithmetic progression. Its quotient genera belong to one of the classes studied by Lee and Nam~\cite{LeeNam2024}. We give a direct calculation suited to Corollary~\ref{cor:quotient-genus}.

\begin{proposition}[Two coprime jumps]\label{prop:two-jump-quotients}
Let $S=\{u,M\}$, where $1\le u<M$ and $\gcd(u,M)=1$. For $q\mid M$, put $N=M/q$, $A=M-u$, and $K=\lfloor A/(2q)\rfloor$. Then
\begin{equation}\label{eq:two-jump-quotient}
\Sigma_S/q=\langle N,M-u,M+u\rangle,
\tag{3.14}\end{equation}
and
\begin{equation}\label{eq:two-jump-quotient-genus}
g(\Sigma_S/q)=\frac{(A-1)(N-1)}2-AK+qK(K+1).
\tag{3.15}\end{equation}
\end{proposition}
\begin{proof}
We have $\Sigma_S=\langle M-u,M,M+u\rangle$. Its elements can be written as $cM+bu$, with $b\in\Z$ and $c\ge|b|$. Since $\gcd(u,q)=1$, divisibility by $q$ forces $b=qj$. After division, the elements are $cN+ju$ with $c\ge q|j|$. For $j\ge0$ this equals $(c-qj)N+j(M+u)$; for $j<0$ it equals $(c-q|j|)N+|j|(M-u)$. This proves~\eqref{eq:two-jump-quotient}.

Index the residue classes modulo $N$ by $ku$, $0\le k<N$, using $\gcd(u,N)=1$. Their Ap\'ery representatives are
\[
w_k=\min\{k(M+u),(N-k)(M-u)\}.
\]
To justify minimality, take an expression for a least element in
the residue class $ku$ modulo $N$. The presence of an $N$-summand
would contradict minimality, since subtracting it preserves the
residue class and stays in the semigroup. Likewise, an expression
using both $M-u$ and $M+u$ cannot be minimal: their sum is
$2M=2qN$, and replacing that pair by $2q$ copies of $N$
would again exhibit a smaller element in the same class after
those copies are deleted. Thus a minimal expression uses only
copies of $M+u$ or only copies of $M-u$.

In the first case, if the expression contains $j\ge0$ copies,
the congruence $ju\equiv ku\pmod N$ gives $j\equiv k\pmod N$;
the least possible choice is $j=k$. In the second case,
$-ju\equiv ku\pmod N$ gives $j\equiv N-k\pmod N$, so the least
choice for $1\le k<N$ is $j=N-k$. For $k=0$, the representative
is $w_0=0$. Any additional complete $N$-cycle of generators
increases the value. This establishes the displayed minimum,
including the case $N=1$.

The first branch is minimal exactly when
$2Mk\le N(M-u)$, or $k\le A/(2q)$.
Set $K=\lfloor A/(2q)\rfloor$; since $0<A<M$,
we have $K<N$. The standard Ap\'ery genus formula~\cite{Rosales}
now gives
\[
\begin{aligned}
g(\Sigma_S/q)
&=\frac{(M+u)K(K+1)+(M-u)(N-K-1)(N-K)}{2N}
-\frac{N-1}{2}\\
&=\frac{(A-1)(N-1)}2-AK+qK(K+1).
\end{aligned}
\]
The same formulas hold for $N=1$, when $K=0$.
\end{proof}

For $u=1$, this simplifies to
\[
g(\Sigma_{\{1,M\}}/q)
=q\left\lfloor\frac{N^2}{4}\right\rfloor
-\left\lfloor\frac N2\right\rfloor.
\]
For example, $S=\{1,6\}$ has defect semigroup $\langle5,6,7\rangle$. The quotient genera for $q=1,2,3,6$ are $6,3,2,0$, respectively. In particular, on even graph orders its eventual dimension is $3D+\lfloor D/2\rfloor-14$. These genus calculations supply explicit constants in the restricted spectral formula.

\section{Exact examples and finite-rank methods}
We next give a few cases where all the filtered dimensions can be written down explicitly. They also provide a useful check on the general results.

For one jump, $S=\{a\}$, introduce the finite fibers
\[
\mathcal F_{a,m}=\{gm:g\mid a,\ \gcd(m,a/g)=1,\ gm>2a\},
\qquad k_{a,m}=|\mathcal F_{a,m}|,\quad m\ge3.
\]
\begin{theorem}[All-degree ranks for one jump]\label{thm:single-jump}
For every $a\ge1$ and $D\ge0$,
\begin{equation}\label{eq:single-rank}
\rho_{\{a\}}^{\glob}(D)=\left\lfloor\frac D2\right\rfloor+1+
\sum_{m=3}^D\min\left\{k_{a,m},1+\left\lfloor\frac{D-m}{2}\right\rfloor\right\}.
\tag{4.1}\end{equation}
The associated-graded global algebra is a square-zero extension of $\Q[t]$, $\deg t=2$, by the modules $\Q[t]/(t^{k_{a,m}})$ placed in weights $m$. In particular,
\[
\rho_{\{1\}}^{\glob}(D)=\left\lfloor\frac D2\right\rfloor+1+\max\{D-2,0\}.
\]
\end{theorem}
\begin{proof}
For one jump, set $\xi_{k,n}=\cos(2\pi ak/n)$ and
$\tau=2(G_2+1)=n$. Let $C_d(X)$ be the Chebyshev polynomial
of the first kind and write $\gamma_d=C_d(0)$.
We denote the M\"obius function on the positive integers by $\mu$.
To avoid mixing an arithmetic function with a symmetric function,
define $\widetilde\tau=2(p_2+1)\in\Lambda_\Q$, so that
$\operatorname{ev}_S(\widetilde\tau)=\tau$. In the
symmetric-function ring put
\[
B_d=\sum_{i=1}^d[X^i]C_d(X)\,p_i
+\gamma_d\widetilde\tau+1-\gamma_d,
\qquad E_m=\frac1m\sum_{d\mid m}\mu(m/d)B_d.
\]
Every term in $B_d$ has weight at most $d$ for $d\ge2$, and $\gamma_1=0$. In particular, $E_m$ has weight at most $m$. The Chebyshev identity $C_d(\cos\theta)=\cos(d\theta)$, followed by root-of-unity orthogonality, gives
\[
\operatorname{ev}_S(B_d)(n)=n\one_{n\mid ad}.
\]
To see the cancellation explicitly, the nonconstant part of $C_d$ sums to
$n\one_{n\mid ad}-1-(n-1)\gamma_d$; the terms added in the definition of $B_d$ cancel the last two contributions.

Let $m_0=n/\gcd(n,a)$. The divisibility $n\mid ad$ is equivalent to $m_0\mid d$. By M\"obius inversion,
\[
e_m(n):=\operatorname{ev}_S(E_m)(n)=\frac nm\one_{m=m_0}
=\frac nm\one_{n\in\mathcal F_{a,m}}.
\]
The coefficient of $p_m$ in $E_m$ is $2^{m-1}/m$, since
$C_m$ has leading coefficient $2^{m-1}$. This proves the
algebraic triangularity needed below. We also claim that
$e_m$, as a function in the filtered image, has \emph{minimum}
weight exactly $m$. It has weight at most $m$ by construction.
On the other hand, $am\in\mathcal F_{a,m}$ for $m\ge3$,
and $e_m(am)=a\ne0$. The first winding length at $am$
is $m$: a signed $a$-walk with nonzero displacement divisible
by $am$ needs at least $m$ steps, and $m$ positive steps
suffice. Every function of weight less than $m$ has, by the
moment formula, a polynomial tail, and its correction vanishes
at $am$, since $am>(m-1)a$. Were $e_m$ of weight less
than $m$, its finite support would force that polynomial tail
to be identically zero; evaluating at $am$ would then give
$e_m(am)=0$, a contradiction.

Now $G_1=-1$ is constant and $G_2$ is affine in $\tau$. The nonzero coefficient of $p_m$ in $E_m$ allows us to solve successively for each $G_m$, $m\ge3$, using $\tau,e_3,\ldots,e_m$ without increasing the weight. Thus these functions generate the global algebra with its filtration.

The fibers $\mathcal F_{a,m}$ partition the admissible levels. We have
\[
e_me_\ell=0\ (m\ne\ell),\qquad
e_m^2=\frac nm e_m,\qquad
\left(\prod_{u\in\mathcal F_{a,m}}(n-u)\right)e_m=0.
\]
The first two relations put every expression in the form of a polynomial in $n$ plus finitely many terms $P_m(n)e_m$. The reduction does not increase weight, because repeated powers of $e_m$ cost $m+2(k-1)\le km$. By the last relation we may divide $P_m$ by a monic polynomial and arrange $\deg P_m<k_{a,m}$, again without increasing weight.

This normal form is unique. First restrict it to an infinite tail: all the $e_m$ vanish there, so the polynomial part is zero. Then restrict to each fiber $\mathcal F_{a,m}$. The function $e_m$ is nonzero on that fiber, while a polynomial of degree less than $k_{a,m}$ cannot vanish at its $k_{a,m}$ distinct points unless it is zero. Hence a filtered basis consists of $\{n^j:2j\le D\}$ together with
\[
\{n^je_m:m+2j\le D,\ 0\le j<k_{a,m}\}.
\]
Counting these vectors proves~\eqref{eq:single-rank}. In the associated-graded algebra, all $e_me_\ell$ have lower-weight representatives, including the case $m=\ell$, and $n^{k_{a,m}}e_m$ also reduces to lower weight. The unique normal form shows that these are exactly the square-zero and truncation relations claimed.
\end{proof}

\begin{theorem}[The first four-regular family]\label{thm:two-jumps}
For $S=\{1,2\}$ and every $D\ge0$,
\begin{equation}\label{eq:two-rank}
\rho_{\{1,2\}}^{\glob}(D)=\left\lfloor\frac D2\right\rfloor+1+
\max\{D-2,0\}+\max\{D-4,0\}.
\tag{4.2}\end{equation}
Moreover,
\[
\gr A_{\{1,2\}}^{\glob}\simeq
\frac{\Q[t,u_3,u_4,\ldots]}
{\bigl(u_hu_k\ (h,k\ge3),\ t^2u_h\ (h\ge3)\bigr)},
\qquad \deg t=2,\ \deg u_h=h.
\]
Its Hilbert series equals $1/(1-z^2)+(z^3+z^5)/(1-z)$.
\end{theorem}
\begin{proof}
For $n\ge5$, any signed $\{1,2\}$-walk with nonzero winding
has length at least $\lceil n/2\rceil$, since the absolute
value of its displacement is at most twice its length.
Conversely, if $n$ is even, $n/2$ steps of $+2$ wind once;
if $n$ is odd, $(n-1)/2$ steps of $+2$ followed by
one step of $+1$ do so. Hence
\[
\nu_{\{1,2\}}(n)=\lceil n/2\rceil\quad(n\ge5),
\qquad B_h=\{2h-1,2h\}\quad(h\ge3).
\]
The two values $Q_h(2h)$ and $Q_h(2h-1)$ are positive,
since the Laurent coefficients of $z^{2h}$ and $z^{2h-1}$
in $L_{\{1,2\}}^h$ are positive. Thus $Q_h,nQ_h$ form
an invertible diagonal-times-Vandermonde pair on $B_h$.
By descending induction,
\[
\{Q_h,nQ_h:3\le h\le H\}
\]
is a basis of all sequences supported on $5\le n\le2H$.

Consequently, in degree at most $D$, the polynomial functions $1,n,\ldots,n^{\lfloor D/2\rfloor}$ together with
$Q_h$ for $3\le h\le D$ and $nQ_h$ for $3\le h\le D-2$
are independent. They also span without an increase in weight.
Indeed, after expanding a polynomial in the $G_h$ into
$n$ and the corrections, a term with a correction has the
form $n^jQ_{h_1}\cdots Q_{h_t}$; put $h=\min h_i$.
Its support lies in $X_h=\{5,\ldots,2h\}$, so triangular
interpolation uses only $Q_k,nQ_k$ for $3\le k\le h$.
The largest required weight is $h+2$. If $t\ge2$,
the original weight is at least $2h\ge h+2$.
If $t=1$ and $j\ge1$, the original weight is at least
$h+2$, and for $j=0$ the term $Q_h$ is already listed.
Thus each reduction respects the original filtration.
Counting gives~\eqref{eq:two-rank}.

In the graded algebra, the products $u_hu_k$ and $t^2u_h$ vanish because their representatives have expressions of smaller weight. The independent filtered basis shows that there are no additional relations.
\end{proof}

\subsection{A finite evaluation-matrix algorithm}
The eventual formulas above have effective but possibly large
thresholds. The following exact finite algorithm applies in
\emph{every} weight, including when the jumps are not coprime
or the modulus does not divide $M$.

\begin{proposition}[Finite evaluation rank, with arbitrary restrictions]
\label{prop:finite-rank}
Let $S$ be any nonempty finite jump set, $D\ge0$, $q\ge1$,
and $0\le a<q$. For this proposition, the notation
$I_{q,a}=\{n>2M:n\equiv a\pmod q\}$ is used for arbitrary $q$.
Set
\[
\begin{aligned}
K&=\max\{2M,DM\},\qquad e=\left\lfloor\frac D2\right\rfloor,\\
\mathcal I_{D;q,a}&=\bigl\{n\in\mathbb Z:2M<n\le K+q(e+1),\\
&\hspace{90pt}n\equiv a\pmod q\bigr\}.
\end{aligned}
\]
For every partition $\lambda=(\lambda_1,\ldots,\lambda_\ell)$
with $|\lambda|=\sum_i\lambda_i\le D$, put
$p_\lambda=\prod_i p_{\lambda_i}$, and set $p_{\varnothing}=1$.
Form the matrix over $\mathbb Q$
\[
\mathcal E_{D;q,a}=
\bigl(\operatorname{ev}_S(p_\lambda)(n)\bigr)_
{n\in\mathcal I_{D;q,a},\,|\lambda|\le D}.
\]
Then
\[
\rank_{\mathbb Q}\mathcal E_{D;q,a}
=\dim_{\mathbb Q}\operatorname{res}_{I_{q,a}}
(A_{S,\le D}^{\glob}).
\]
In particular, $q=1$ computes $\rho_S^{\glob}(D)$ exactly.
\end{proposition}
\begin{proof}
The columns of the matrix are evaluations of a spanning family
of $A_{S,\le D}^{\glob}$, so its rank cannot exceed the
dimension of the restriction image. To prove the reverse
inequality, it suffices to show that restriction from that
image to the finite set $\mathcal I_{D;q,a}$ is injective.

Let $f$ be an image of weight at most $D$ vanishing on this set.
Every moment $G_h$ occurring in a weight-$D$ monomial has
$h\le D$. For $n>K\ge DM$, all these moments are affine
in $n$ by Lemma~\ref{lem:global-bridge}. As in the proof of
Proposition~\ref{prop:glob-image}, the restriction of $f$ to
$n>K$ therefore agrees with a polynomial $P(n)$ of degree
at most $e$. There are at least $e+1$ integers in the
specified residue class with
$K<n\le K+q(e+1)$, and all lie in
$\mathcal I_{D;q,a}$. Hence $P$ has at least $e+1$
distinct zeros, forcing $P=0$.
For every remaining admissible $n\le K$ in that residue class,
$n\in\mathcal I_{D;q,a}$ and $f(n)=0$ by assumption.
Thus $f$ vanishes on the whole of $I_{q,a}$, proving
injectivity and the rank formula.
\end{proof}

\section{Discussion and further questions}
The filtration of polynomial spectral moments measures the cost of
detecting winding walks in a family of circulant graphs. Its eventual
dimension records the genus of $\Sigma_S$, while restrictions of
the graph order to congruence classes recover individual gap counts.
When all classes modulo $M$ are retained, one stable-degree list of
dimensions gives the Ap\'ery set and reconstructs the semigroup.

Two natural directions remain. First, one may seek conditions under
which a smaller family of congruence restrictions still distinguishes
the defect semigroups arising from jump sets. Second, it would be useful
to characterize intrinsically which numerical semigroups can arise in
the form $\langle M,M-s,M+s:s\in S\rangle$ and which spectral data
are necessary for their reconstruction. These are questions about the
information carried by the filtered graph invariants, rather than
about quotient genera in isolation.

\section*{Statements and Declarations}

\noindent\textbf{Funding.}
The authors received no funding for this research.

\medskip

\noindent\textbf{Competing interests.}
The authors declare that they have no relevant
financial or non-financial interests to disclose.

\medskip

\noindent\textbf{AI-assisted preparation.}
OpenAI ChatGPT (GPT-6) was used as an aid in drafting and revision of arguments,
and in LaTeX and language editing.
AI is not an author.
The named authors take full responsibility for
checking the mathematical arguments, verifying
the references, and approving the submitted
version of this manuscript.
\end{document}